\documentclass[12pt,letterpaper]{article}

\usepackage{latexsym} 
\usepackage{array}
\usepackage{enumerate}
\usepackage{theorem}
\newtheorem{lemma}{Lemma}

\newtheorem{theorem}{Theorem}

\begin{document}

\title{BOSONIC KERNELS}
\author{P. L. Robinson}
\date{}
\maketitle 

\begin{abstract}
We show that the bosonic Fock representation of a complex Hilbert space admits a purely algebraic kernel calculus; as an illustration, we use it to reproduce the standard integral kernel formulae for metaplectic operators within the complex-wave representation.  
\end{abstract}

\section*{Introduction}

One of the most elegant models of the bosonic Fock representation of a complex Hilbert space $V$ is the complex-wave model, traditionally due to Bargmann \cite{Barg1} in finite dimensions and to Segal \cite{Segal} \cite{Segalcw} in infinite dimensions. In this formulation, the Fock space $\mathbf{H}$ comprises all square-integrable antiholomorphic functions $V \longrightarrow \mathbf{C}$, square-integrability being defined relative to a coherent family of Gaussian probability measures $\mu$ on finite-dimensional complex subspaces and antiholomorphicity being suitably tamed by the same subspaces. A particularly attractive feature of this model is the fact that each bounded linear operator $T$ on $\mathbf{H}$ is an integral operator with kernel $t: V \times V \longrightarrow \mathbf{C}$ that is antiholomorphic-holomorphic and satisfies 
\nonumber 
\[
f \in \mathbf{H}; z \in V \Longrightarrow Tf(z) = \int_{V} t(z,w)f(w)d\mu(w). 
\]

\medbreak 

This has been used to great effect in (for instance) making explicit the associated metaplectic representation \cite{Shale} of the restricted symplectic group, comprising all real symplectic automorphisms $g$ of $V$ for which the commutator $[g, i]$ is of Hilbert-Schmidt class: each such $g$ induces on $\mathbf{H}$ an implementing metaplectic operator $U$ that is unitary up to scale and has integral kernel $u$ given by 
\nonumber 
\[
z,w \in V \Longrightarrow u(z,w) = \exp \frac{1}{2}\Bigl \{2\langle C_{g}^{-1} z \vert w \rangle + \langle z \vert Z_{g^{-1}} z \rangle + \langle Z_g w \vert w \rangle \Bigr \} 
\]
where $C_g = \frac{1}{2} (g - i g i)$ is the (automatically invertible) complex-linear part of $g$ and $Z_g = C_{g}^{-1} (g - C_g)$ is antilinear: for example, see \cite{Barg2} \cite{Berezin} \cite{Foll} \cite{Neretin} \cite{RobRaw} \cite{Vergne} . More recently, it has been shown that the foregoing kernel defines in $\mathbf{H}$ a sesquilinear form that weakly implements $g$ quite generally, without the Hilbert-Schmidt condition: see \cite{PPSZ} and \cite{Rob}.

\medbreak 

Our chief aim in this paper is to demonstrate that, when viewed correctly, the standard bosonic Fock representation of $V$ accommodates a purely algebraic kernel calculus with remarkably pleasant features. As a particular illustration, we rather easily recover precisely the kernels of the metaplectic operators, in both the weak form and the strong. At the same time, we take the opportunity to treat in like fashion the implementation of antisymplectics. 

\medbreak 

In retrospect, the bosonic kernel calculus presented in this paper may be seen to have origins in the work of Berezin \cite{Berezin}. There, the ideas are somewhat veiled by (among other things) the assumption that $V$ be a space of square-integrable functions and the related involvement of coordinates. Here, matters are reduced to the invariant essentials, promoting clarity and further development. 

\section*{Fock space}

Our aim in this section is to present as background well-known bosonic Fock-lore in a form that is specially tailored to the subsequent development of bosonic kernels and their applications. Details may be found in the standard mathematical references on quantum field theory; the precise approach adopted here originally appeared in \cite{Rob} to which we refer for proofs. 
\medbreak 

Let $V$ be a complex Hilbert space with $\langle \cdot \vert \cdot \rangle$ as its inner product. Denote by $SV = \bigoplus_{n = 0}^{\infty } S^n V$ its (graded) symmetric algebra and by $SV'$ the full antidual, comprising all antilinear functionals $SV \longrightarrow \mathbf{C}$. Note that $SV$ carries a natural cocommutative coproduct: the diagonal map $V \longrightarrow V \oplus V$ induces an algebra homomorphism $SV \longrightarrow S(V \oplus V)$ which when followed by the canonical isomorphism $S(V \oplus V) \equiv SV \otimes SV$ yields $\Delta : SV \longrightarrow SV \otimes SV$. In consequence, the full antidual $SV'$ carries a natural commutative product defined by the prescription that if $\Phi, \Psi \in SV'$ then 
\nonumber 
\[
\theta \in SV \Longrightarrow [\Phi \Psi](\theta) = [\Phi \otimes \Psi](\Delta \theta).
\]
In addition $SV$ carries a canonical inner product, determined by perpendicularity of the homogeneous summands $(S^n V : n \in \mathbf{N})$ along with the requirement that if $n \in \mathbf{N}$ and $x_1, \dots , x_n, y_1, \dots , y_n \in V$ then 
\nonumber 
\[
\langle x_1 \cdots x_n \vert y_1 \cdots y_n \rangle = \sum_{p} \prod_{j = 1}^{n} \langle x_j \vert y_{p(j)} \rangle 
\]
where $p$ runs over all permutations of $\{ 1, \dots , n \}$. These structures are related by the fact that a canonical algebra embedding is defined by the rule 
\nonumber 
\[
SV \longrightarrow SV' : \phi \longmapsto \langle \cdot \vert \phi \rangle.
\]
Finally, bosonic Fock space $S[V] = \bigoplus_{n = 0}^{\infty } S^n [V]$ itself may be defined either (abstractly) as the Hilbert space completion of $SV$ or (concretely) as the subspace of $SV'$ comprising all bounded antifunctionals. 

\medbreak 

The bosonic Fock representation of $V$ is fashioned from creators and annihilators, defined as follows. Let $v \in V$. The creator $c(v) : SV \longrightarrow SV$ is the operator of multiplication by $v$: thus, 
\nonumber 
\[
\phi \in SV \Longrightarrow c(v) \phi = v \phi.
\] 
The annihilator $a(v) : SV \longrightarrow SV$ is the unique linear derivation that kills $\mathbf{1} \in \mathbf{C} = S^0 V$ and sends $w \in V = S^1 V$ to $\langle v \vert w \rangle$: thus, if $v_1, \dots , v_n \in V$ then 
\nonumber 
\[
a(v) (v_1 \cdots v_n) = \sum_{j = 1}^{n} \langle v \vert v_j \rangle v_1 \cdots \widehat{v_j} \cdots v_n
\]
where the circumflex signifies omission as usual. These operators are mutually adjoint: if $\phi, \psi \in SV$ then 
\nonumber 
\[
\langle \psi \vert c(v) \phi \rangle = \langle a(v) \psi \vert \phi \rangle
\]
as may be conveniently verified by multilinearity and evaluation on decomposables. Accordingly, they extend to $SV'$ by antiduality: thus, if $\Phi \in SV'$ and $\psi \in SV$ then 
\nonumber 
\[
[c(v) \Phi] (\psi) = \Phi [a(v) \psi] 
\] 
and 
\nonumber 
\[
[a(v) \Phi] (\psi) = \Phi [c(v) \psi];
\]
note that $c(v)$ continues to be a multiplication operator and $a(v)$ continues to be a linear derivation. We remark that if $v$ is nonzero then $c(v)$ is injective; we remark also that the intersection of the kernels of all annihilators is the complex line spanned by the functional $\langle \cdot \vert \mathbf{1} \rangle \in SV'$. Finally, the corresponding mutually adjoint creators and annihilators in Fock space $S[V] \subset SV'$ are defined by restriction to the (coincident) natural domains 
\nonumber 
\[
\mathcal{D}[c(v)] = \{ \Phi \in S[V] : c(v) \Phi \in S[V] \}
\]
\nonumber 
\[
\mathcal{D}[a(v)] = \{ \Phi \in S[V] : a(v) \Phi \in S[V] \}.
\]

\medbreak 

Creators and annihilators satisfy the canonical commutation relations: if $x, y \in V$ then 
\nonumber 
\[
[a(x) , a(y)] = [c(x) , c(y)] = 0
\]
\nonumber 
\[
[a(x) , c(y)] = \langle x \vert y \rangle I.
\]
On $SV$ these relations hold by direct calculation: for the mixed commutator, if $\psi \in SV$ then 
\nonumber 
\[
a(x)c(y) \psi = [a(x)y] \psi + y [a(x) \psi] = \langle x \vert y \rangle \psi + c(y)a(x) \psi.
\]
On $SV'$ precisely the same relations continue to hold by antiduality. The situation on Fock space $S[V]$ itself is more delicate, domain considerations forcing the replacement of equalities by inclusions; for instance, if $x,y \in V$ then $[a(x) , c(y)] \subset \langle x \vert y \rangle I$. 

\medbreak 

The various structures considered thus far admit complexification. For this purpose, it helps to regard $V$ as a real vector space with complex structure $J = i \cdot$ and introduce the real symplectic form $\Omega = \mathrm{Im} \langle \cdot \vert \cdot \rangle$ so that 
\nonumber 
\[
x,y \in V \Longrightarrow  \langle x \vert y \rangle = \Omega (x,J y) + i \Omega (x,y).
\]
Now, let $V_{\mathbf{C}} = \mathbf{C} \otimes V$ be the complexification of $V$ with overbar $\bar{\cdot}$ as the canonical conjugation. The complexifications of $J$ and $\Omega$ furnish $V_{\mathbf{C}}$ with a (complete) inner product according to the prescription 
\nonumber 
\[
x,y \in V_{\mathbf{C}}\Longrightarrow \langle x \vert y \rangle = \Omega_{\mathbf{C}} (\bar{x}, J_{\mathbf{C}} y).
\]
The eigendecomposition $V_{\mathbf{C}} = V_{+} \oplus V_{-}$ with $J_{\mathbf{C}} \vert V_{\pm} = \pm  i I$ is orthogonal, with induced isometry 
\nonumber 
\[
V \longrightarrow V_{+} : v \longmapsto v_{+} = \frac{1}{\sqrt{2}}(v - i J v)
\]
and induced antiisometry 
\nonumber 
\[
V \longrightarrow V_{-} : v \longmapsto v_{-} = \frac{1}{\sqrt{2}}(v + i J v)
\]
so that if $x,y \in V$ then $\langle x_{+} \vert y_{+} \rangle = \langle x \vert y \rangle $ and $\langle x_{-} \vert y_{-} \rangle = \langle y \vert x \rangle $. 

\medbreak 

The various notions of Fock space also complexify. Thus, $SV_{\mathbf{C}}$ acquires a canonical inner product and thereby embeds in its full antidual $SV_{\mathbf{C}}'$ while the Fock space $S[V_{\mathbf{C}}]$ sits between: $SV_{\mathbf{C}} \subset S[V_{\mathbf{C}}] \subset SV_{\mathbf{C}}'$. Regarding symmetric algebras, the isometry $V \longrightarrow V_{+} \subset V_{\mathbf{C}}$ extends to an isometric algebra isomorphism $SV \stackrel{+}{\longrightarrow} S(V_{+}) \subset SV_{\mathbf{C}}$ while the antiisometry $V \longrightarrow V_{-} \subset V_{\mathbf{C}}$ extends to an antiisometric algebra (anti-)isomorphism $SV \stackrel{-}{\longrightarrow} S(V_{-}) \subset SV_{\mathbf{C}}$; in this connexion, note that the direct sum decomposition $V_{\mathbf{C}} = V_{+} \oplus V_{-}$ engenders a tensor product decomposition $SV_{\mathbf{C}} \equiv S(V_{+}) \otimes S(V_{-})$. Regarding full antiduals, the orthogonal projection $V_{\mathbf{C}} \longrightarrow V_{+}$ induces a homomorphism $SV_{\mathbf{C}} \longrightarrow S(V_{+}) \longrightarrow SV$ which passes by antiduality to a homomorphism 
\nonumber 
\[
SV' \longrightarrow SV_{\mathbf{C}}' : \Phi \longmapsto \Phi_{+} 
\]
satisfying the compatibility condition 
\nonumber 
\[
\Phi \in SV' ; \psi \in SV \Longrightarrow \Phi_{+} (\psi_{+} ) = \Phi (\psi); 
\]
likewise, there is a canonical antilinear homomorphism 
\nonumber 
\[
SV' \longrightarrow SV_{\mathbf{C}}' : \Phi \longmapsto \Phi_{-} 
\]
satisfying  
\nonumber 
\[
\Phi \in SV' ; \psi \in SV \Longrightarrow \Phi_{-} (\psi_{-}) = \overline{\Phi (\psi)}. 
\]

\medbreak 

Gaussians play a significant role in the development of Fock space. Denote by $S^2 V'$ the space of all quadratics on $V$: elements $\zeta$ of $SV'$ that are concentrated in second degree, in the sense that if $n \neq 2$ then $\zeta \vert S^n V = 0$. Denote by $\Sigma V '$ the space of all (necessarily) antilinear maps $Z : V \longrightarrow V' $ that are symmetric in the sense that if $x,y \in V$ then $Z x (y) = Z y (x)$. A canonical linear isomorphism 
\nonumber 
\[
S^2 V' \longrightarrow \Sigma V ' : \zeta \longmapsto Z
\]
is defined by the rule 
\nonumber 
\[
x,y \in V \Longrightarrow Z x (y) = \zeta (x y)
\]
whence 
\nonumber 
\[
v \in V \Longrightarrow a(v) \zeta = Z v.
\]
Under this isomorphism, $S^2 [V]$ corresponds to the Hilbert space $\Sigma[V]$ comprising all Hilbert-Schmidt antilinear operators $Z : V \longrightarrow V$ that are symmetric in the sense that if $x,y \in V$ then $\langle y \vert Z x \rangle = \langle x \vert Z y \rangle$; moreover, if $Z \in \Sigma [V]$ and $\zeta \in S^2 [V]$ correspond, then $\Vert Z \Vert_{HS} = \sqrt{2}\Vert \zeta \Vert $. 

\medbreak 

Now, let $Z \in \Sigma V'$ and $\zeta \in S^2 V'$ correspond: the associated Gaussian $ e^{Z} = \exp\zeta$ is defined by the power series 
\nonumber 
\[
e^{Z} = \sum_{n = 0}^{\infty} \frac{1}{n!} \zeta^n \in SV'
\]
which always converges pointwise on $SV$. Note that annihilators act essentially as creators on a Gaussian: if $Z \in \Sigma V'$ then 
\nonumber 
\[
v \in V \Longrightarrow a(v) e^Z = (Z v) e^Z. 
\]
We remark that $e^Z$ lies in Fock space $S[V]$ precisely when $Z$ lies in $\Sigma[V]$ and has operator norm strictly less than unity: in this case, $e^Z$ actually lies in the domain of each polynomial in creators and annihilators; moreover, if $\mathrm{Det}$ denotes Fredholm determinant then 
\nonumber 
\[
\langle e^Z \vert e^Z \rangle = \mathrm{Det}^{\frac{1}{2}}(I - Z^2)^{-1}.
\]

\medbreak 

In particular, the extended algebra $SV_{\mathbf{C}}'$ contains a preferred Gaussian: the canonical conjugation on $V_{\mathbf{C}}$ is symmetric antilinear, for if $x,y \in V_{\mathbf{C}}$ then $\langle y \vert \overline{x} \rangle = \langle x \vert \overline{y} \rangle$; thus it corresponds to a canonical quadratic $\zeta_V \in S^2 V_{\mathbf{C}}'$ according to 
\nonumber 
\[
x,y \in V_{\mathbf{C}} \Longrightarrow \zeta_V (xy) = \langle y \vert \overline{x} \rangle. 
\]
The preferred Gaussian $\exp \zeta_V = e^{Z_V} \in SV_{\mathbf{C}}'$ features in our development of bosonic kernels. 

\medbreak 

At this point, it is convenient to insert a useful property of the canonical inner product on the symmetric algebra. In fact, let $[ \cdot \vert \cdot]$ be a Hermitian form on $SV$ with the property that 
\nonumber 
\[
v \in V;  \phi, \psi \in SV \Longrightarrow [c(v) \psi \vert \phi] = [\psi \vert a(v) \phi].
\]
We claim that this property characterizes the canonical inner product, up to normalization: $[ \cdot \vert \cdot]$ is proportional to $\langle \cdot \vert \cdot \rangle$. To see this, note first that if $v_1, \dots ,v_n \in V$ and $\phi \in S^n V$ then 
\nonumber 
\[
a(v_n) \cdots a(v_1) \phi = \langle v_1 \cdots v_n \vert \phi \rangle  
\]
as follows readily by applying first to the case $\phi = v^n$ with $v \in V$ and noting that $S^n V$ is spanned by such powers, whence by hypothesis, 
\nonumber 
\[
[v_1 \cdots v_n \vert \phi] = [\mathbf{1} \vert a(v_n) \cdots a(v_1) \phi] = [\mathbf{1} \vert \mathbf{1} ] \langle v_1 \cdots v_n \vert \phi \rangle; 
\]
likewise, if $\phi \in S^m V$ with $m < n$ then $a(v_n) \cdots a(v_1) \phi = 0$ whence it follows that $[v_1 \cdots v_n \vert \phi] = 0$. Finally, the Hermitian nature of $[ \cdot \vert \cdot]$ (or a symmetric argument) completes the justification of our claim. 

\medbreak 

\begin{lemma} \label{innerproduct}
If the Hermitian form $[ \cdot \vert \cdot]$ on $SV$ satisfies the property 
\nonumber 
\[
v \in V;  \phi, \psi \in SV \Longrightarrow [c(v) \psi \vert \phi] = [\psi \vert a(v) \phi]
\]
then 
\nonumber 
\[
[ \cdot \vert \cdot] = [\mathbf{1} \vert \mathbf{1} ]\langle \cdot \vert \cdot \rangle. 
\]
\end{lemma} 

\begin{flushright}
$\Box$
\end{flushright} 

\section*{Kernel calculus}

Our aim in this section is to present an analogue of the Schwarz kernel theorem from generalized function theory, whereby suitable linear maps are represented by ``integral'' kernels. 

\medbreak 

To begin, let $u \in SV_{\mathbf{C}}'$. If $\phi, \psi \in SV$ then $\psi_{+} \phi_{-} \in S(V_{+}) \otimes S(V_{-}) \equiv SV_{\mathbf{C}}$ whence $u ( \psi_{+} \phi_{-} ) \in \mathbf{C}$ is defined. In this way, to $u \in SV_{\mathbf{C}}'$ is associated a linear map $U: SV \longrightarrow SV'$ given by the rule that if $\phi, \psi \in SV$ then $[U \phi] (\psi) = u ( \psi_{+} \phi_{-} )$. In the opposite direction, if $U: SV \longrightarrow SV'$ is linear then the map 
\nonumber 
\[
S(V_{+}) \times S(V_{-}) \longrightarrow \mathbf{C} : (\psi_{+} , \phi_{-}) \longmapsto [U \phi] (\psi)
\]
is bi-antilinear, hence descends to an antilinear functional $u \in SV_{\mathbf{C}}'$ on  $S(V_{+}) \otimes S(V_{-}) \equiv SV_{\mathbf{C}}$. The correspondence $u \longleftrightarrow U$ being plainly bijective and linear in each direction, we have established the following kernel theorem. 

\medbreak 

\begin{theorem} \label{kernel}
A natural linear isomorphism 
\nonumber 
\[
SV_{\mathbf{C}}' \longrightarrow \mathrm{Hom}(SV, SV') : u \longmapsto U
\]
is defined by the prescription 
\nonumber 
\[
\phi, \psi \in SV \Longrightarrow [U \phi] (\psi) = u (\psi_{+} \phi_{-}).
\]
\end{theorem} 

\begin{flushright}
$\Box$
\end{flushright} 

\medbreak 

We refer to $u$ as the kernel of $U$. 

\medbreak 

The canonical conjugation on $V_{\mathbf{C}}$ extends to an involution on the algebra $SV_{\mathbf{C}}$ which passes to $SV_{\mathbf{C}}'$ by antiduality: 
\nonumber 
\[
u \in SV_{\mathbf{C}}' ; \theta \in SV_{\mathbf{C}} \Longrightarrow u^* (\theta) = \overline{u (\theta^*)}.
\]
Each linear map $U: SV \longrightarrow SV'$ has a canonical adjoint $U^* : SV \longrightarrow SV'$ defined by 
\nonumber 
\[
\phi, \psi \in SV \Longrightarrow [U^* \phi] (\psi) = \overline{[U \psi] (\phi)}. 
\]
Now, let $U: SV \longrightarrow SV'$ have $u \in SV_{\mathbf{C}}'$ as kernel: if $\phi, \psi \in SV$ then 
\nonumber 
\[
[U^* \phi](\psi) = \overline{[U\psi](\phi)} = \overline{u ( \phi_{+} \psi_{-} )} = \overline{u ( (\psi_{+} \phi_{-})^* )} = u^*( \psi_{+} \phi_{-} )
\]
whence $U^*$ has $u^*$ as kernel. This proves the following. 

\medbreak 

\begin{theorem} \label{adjoint}
If $U \in \mathrm{Hom}(SV, SV')$ has kernel $u \in SV_{\mathbf{C}}'$ then $U^*$ has kernel $u^*$. 
\end{theorem} 

\begin{flushright}
$\Box$
\end{flushright} 

\medbreak 

Let us now consider some elementary examples, starting with the kernel of the identity operator $I: SV \longrightarrow SV \subset SV'$. For this, recall that the canonical conjugation in $V_{\mathbf{C}}$ prefers the quadratic $\zeta_V \in S^2 V_{\mathbf{C}}'$ satisfying 
\nonumber 
\[
w \in V_{\mathbf{C}} \Longrightarrow a(w) \zeta_V = \overline{w}.
\] 
Now, if $x,y \in V_{\mathbf{C}}$ then 
\nonumber 
\[
\zeta_V (x_{+} y_{-} ) = a(x_{+}) \zeta_V (y_{-}) = x_{-} (y_{-}) = \langle y_{-} \vert x_{-} \rangle = \langle x \vert y \rangle.
\]
More generally, we claim that if $n \in \mathbf{N}$ and $\xi, \eta \in S^n V$ then 
\nonumber 
\[
\zeta_{V}^n (\xi_{+} \eta_{-}) = n! \: \langle \xi \vert \eta \rangle.
\]
To see this, let $x,y \in V$ and proceed inductively: 
\begin{eqnarray*}
\zeta_{V}^{n + 1} ( (x^{n + 1})_{+} (y^{n + 1})_{-} ) & = & a(x_{+}) \zeta_{V}^{n + 1} (x_{+}^{n} y_{-}^{n + 1}) \\ & = & (n + 1) c(x_-) \zeta_{V}^n (x_{+}^{n} y_{-}^{n + 1}) \\ & = & (n + 1) \zeta_{V}^n(a(x_-) \{x_{+}^{n} y_{-}^{n + 1} \}) \\ & = & (n + 1) \zeta_{V}^n ((n + 1) \langle x_- \vert y_- \rangle x_{+}^{n} y_{-}^{n}) \\ & = & (n + 1)^2 \langle x \vert y \rangle \zeta_{V}^n (x_{+}^{n} y_{-}^{n}) \\ & = & (n + 1)^2 \langle x \vert y \rangle n! \langle x^n \vert y^n \rangle \\ & = & (n + 1)! \langle x^{n + 1} \vert y^{n + 1} \rangle 
\end{eqnarray*}
whereupon the claim follows as $S^n V$ is spanned by $n$-fold powers. Summing up, we have established that the kernel of the identity operator on $SV$ is none other than the preferred Gaussian in $SV_{\mathbf{C}}'$. 

\medbreak 

\begin{theorem} \label{identity}
The kernel of the identity operator $I : SV \longrightarrow SV \subset SV'$ is the Gaussian $e^{Z_V} \in SV_{\mathbf{C}}'$: 
\nonumber 
\[
\phi, \psi \in SV \Longrightarrow \langle \psi \vert \phi \rangle = e^{Z_V} (\psi_{+} \phi_{-}).
\]
\end{theorem} 

\begin{flushright}
$\Box$
\end{flushright} 

\medbreak 

We remark that this result offers an alternative (perhaps more elegant) route to defining the canonical inner product on the symmetric algebra. 

\medbreak 

Now, let $U: SV \longrightarrow SV'$ have kernel $u \in SV_{\mathbf{C}}'$ and let $v \in V$.  If $\phi, \psi \in SV$ then 
\begin{eqnarray*}
[(c(v) U) \phi] (\psi) & = & U \phi (a(v) \psi) \\ & = & u ( (a(v) \psi)_{+} \phi_{-}) \\ & = & u ((a(v_{+}) \psi_{+}) \phi_{-}) \\ & = & u(a(v_{+}) (\psi_{+} \phi_{-})) \\ & = & [c(v_{+}) u](\psi_{+} \phi_{-})
\end{eqnarray*}
whence it follows that the operator $c(v) U$ has kernel $c(v_{+}) u$. In precisely similar fashion we may determine the kernels of the remaining three products of $U$ with creators and annihilators, as follows. 

\medbreak 

\begin{theorem} \label{acU}
If $U: SV \longrightarrow SV'$ has kernel $u \in SV_{\mathbf{C}}'$ and if $v \in V$ then the following operators and kernels correspond: 
\nonumber 
\[
c(v) U \longleftrightarrow c(v_+) u
\]
\nonumber 
\[
a(v) U \longleftrightarrow a(v_+) u 
\]
\nonumber 
\[
U c(v) \longleftrightarrow a(v_-) u
\]
\nonumber 
\[
U a(v) \longleftrightarrow c(v_-) u.
\]
\end{theorem} 

\begin{flushright}
$\Box$
\end{flushright} 

\medbreak 

In particular, as the identity operator $I$ has the preferred Gaussian $e^{Z_V}$ as its kernel, if $v \in V$ then the creator $c(v)$ has kernel $(v_+) e^{Z_V}$ and the annihilator $a(v)$ has kernel $(v_-) e^{Z_V}$. Note that the mutually adjoint nature of creators and annihilators is reflected in the mutually adjoint nature of their kernels. 

\medbreak 

Recall that the number operator $N$ is defined initially on $SV$ by 
\nonumber 
\[
n \in \mathbf{N} \Longrightarrow N \vert S^n V = n I
\]
and on $SV'$ by antiduality, thus 
\nonumber 
\[
\Phi \in SV' ; \psi \in SV \Longrightarrow [N \Phi] (\psi) = \Phi [N \psi]. 
\]
If $n > 0$ then 
\nonumber 
\[
\zeta_V  \frac{\zeta_{V}^{n - 1}}{(n - 1)!} = \frac{\zeta_{V}^n}{(n - 1)!} = n \frac{\zeta_{V}^n}{n!}
\]
whence if also $\phi, \psi \in S^n V$ then 
\nonumber 
\[
\zeta_V  e^{Z_V} (\psi_+ \phi_-) = n \frac{\zeta_{V}^n}{n!}(\psi_+ \phi_-) = n \langle \psi \vert \phi \rangle = \langle \psi \vert N\phi \rangle
\]
which justifies the following result. 

\medbreak 

\begin{theorem} \label{number}
The number operator $N$ has kernel $\zeta_V e^{Z_V}$. 
\end{theorem}

\begin{flushright}
$\Box$
\end{flushright} 

\medbreak 

As a further elementary example, consider rank-one operators. Thus, let $\Phi, \Psi \in SV'$: if also $\phi, \psi \in SV$ then 
\nonumber 
\[
(\Psi_+ \Phi_-) (\psi_+ \phi_-) = \Psi_+(\psi_+) \overline{\Phi_-(\phi_-)} = \Psi(\psi)\overline{\Phi(\phi)}. 
\]
Accordingly, $\Psi_+ \Phi_- \in SV_{\mathbf{C}}'$ is the kernel of the linear operator 
\nonumber 
\[
SV \longrightarrow SV' : \phi \longmapsto \overline{\Phi(\phi)}\Psi.
\]

\medbreak 

The kernel calculus developed thus far has matched linear maps $SV \longrightarrow SV'$ with kernels in $SV_{\mathbf{C}}'$. This kernel calculus may be modified to deal with other classes of maps; we shall be content to indicate a modification appropriate to antilinear maps $SV \longrightarrow SV'$. 

\medbreak 

A suitable receptacle for the kernels of antilinear maps $SV \longrightarrow SV'$ is the full antidual $S(V \oplus V)'$ where $V \oplus V$ denotes the direct sum of two copies of $V$ as a complex Hilbert space. Note that the first and second inclusions $V \longrightarrow V \oplus V$ induce algebra embeddings 
\nonumber 
\[
SV' \longrightarrow S(V \oplus V)' : \Phi \longmapsto \Phi_1 
\]
\nonumber 
\[ 
SV' \longrightarrow S(V \oplus V)' : \Psi \longmapsto \Psi_2 
\] 
which restrict to corresponding maps $SV \longrightarrow S(V \oplus V)$. In these terms, an argument along the same lines as that for Theorem \ref{kernel} has the following result. 

\medbreak 

\begin{theorem} \label{antikernel}
The correspondence $u \leftrightarrow U$ set up by the rule  
\nonumber 
\[ 
\phi, \psi \in SV \Longrightarrow [U \phi] (\psi) = u(\phi_1 \psi_2) 
\] 
is a canonical linear isomorphism between $S(V \oplus V)'$ and the space of all antilinear maps $SV \longrightarrow SV'$. 
\end{theorem} 

\begin{flushright}
$\Box$
\end{flushright} 

\medbreak 

We leave the detailed elaboration of this calculus as an exercise, merely remarking that if the antilinear $U:SV \longrightarrow SV'$ has kernel $u \in S(V \oplus V)'$ and if $v \in V$ then (for example) $U c(v)$ has kernel $a(v_1) u$ while $c(v) U$ has kernel $c(v_2) u$.  

\section*{Metaplectic operators}

Our primary aim in this section is to present a new, explicit analysis of the (weak) implementation problem for symplectic automorphisms in the Fock representation, as an extended illustration of the use of bosonic kernels. 

\medbreak 

The symplectic group $\mathrm{Sp}(V)$ comprises all real-linear automorphisms $g$ of $V$ that are symplectic in the sense 
\nonumber 
\[
x,y \in V \Longrightarrow \Omega(g x,g y) = \Omega(x,y).
\]
When $g \in \mathrm{Sp}(V)$ we write $C_g = \frac{1}{2}(g - JgJ)$ for its complex-linear part and $A_g =  \frac{1}{2}(g + JgJ)$ for its antilinear part; we remark that $C_g$ is invertible and that the antilinear operator $Z_g = C_g^{-1} A_g = -A_{g^{-1}} C_{g^{-1}}^{-1}$ is symmetric with operator norm strictly less than unity. 

\medbreak 

The bosonic Fock representation of $V$ is in fact constructed as follows. We define 
\nonumber 
\[
\pi : V \longrightarrow \mathrm{End}\:  SV
\]
by declaring that if $v \in V$ then 
\nonumber 
\[
\pi (v) = \frac{1}{\sqrt{2}} \{ c(v) + a(v) \}
\]
and extend to 
\nonumber 
\[
\pi : V \longrightarrow \mathrm{End}\:  SV'
\]
by antiduality, so that if also $\Phi \in SV'$ and $\psi \in SV$ then 
\nonumber 
\[
[\pi (v) \Phi] (\psi) = \Phi [\pi(v) \psi].
\]
On account of the canonical commutation relations satisfied by creators and annihilators, $\pi$ satisfies the canonical commutation relations in Heisenberg form: the relations 
\nonumber 
\[
x,y \in V \Longrightarrow [\pi(x),\pi(y)] = i \Omega(x,y) I
\]
hold on both $SV$ and $SV'$ precisely as stated. The bosonic Fock representation proper is by operators in Fock space $S[V] \subset SV'$: thus, for $v \in V$ we define $\pi(v)$ by restriction to the natural domain 
\nonumber 
\[
\mathcal{D}[\pi(v)] = \{ \Phi \in S[V] : \pi(v) \Phi \in S[V] \}
\]
on which the resulting operator is selfadjoint. For these operators, domain considerations weaken the Heisenberg form of the canonical commutation relations: if $x,y \in V$ then $[\pi(x),\pi(y)] \subset i \Omega(x,y) I$. 

\medbreak 

Now, if $g \in \mathrm{Sp}(V)$ then $\pi \circ g$ and $\pi$ satisfy the same Heisenberg canonical commutation relations, whence it is reasonable to ask whether they are equivalent in any sense. The traditional sense is in terms of unitary intertwining operators on Fock space $S[V]$ and leads to the famous metaplectic group via the celebrated Shale theorem. We shall consider equivalence in a weak (though eminently reasonable) sense. Specifically, we ask whether there exists a (nonzero) linear map $U : SV \longrightarrow SV'$ that intertwines $\pi$ on $SV$ and $\pi \circ g$ on $SV'$ in the sense that if $v \in V$ then 
\nonumber 
\[
 U \pi(v) = \pi(g v) U
\]
or equivalently 
\nonumber 
\[
U c(v) = \{ c(C_g v) + a(A_g v) \} U 
\]
\nonumber 
\[
U a(v) = \{ c(A_g v) + a(C_g v) \} U. 
\]
As we shall see, such a weak metaplectic operator always exists and is unique up to scalar multiples; our proof is streamlined by the use of bosonic kernels. 

\medbreak 

Thus, let $g \in \mathrm{Sp}(V)$ and let $U : SV \longrightarrow SV'$ be a corresponding weak metaplectic operator with kernel $u \in SV_{\mathbf{C}}'$. We begin by expressing the foregoing relations between $U$, creators and annihilators on $SV \subset SV'$ in terms of $u$, creators and annihilators on $SV_{\mathbf{C}} \subset SV_{\mathbf{C}}'$: by Theorem \ref{acU}, these relations respectively become 
\nonumber 
\[
a(v_-) u = c(C_g v)_+ u + a(A_gv)_+ u 
\] 
\nonumber 
\[
c(v_-) u = c(A_g v)_+ u + a(C_g v)_+ u. 
\]

\medbreak 

Let us seek as solution $u$ a Gaussian $e^Z$: in view of the fact 
\nonumber 
\[
w \in V_{\mathbf{C}} \Longrightarrow a(w) e^Z = (Z w) e^Z 
\]
the conditions on $u = e^Z$ respectively become that if $v \in V$ then 
\nonumber 
\[
Z(v_-) = (C_g v)_+ + Z(A_g v)_+ 
\]
\nonumber 
\[
v_- = (A_g v)_+ + Z(C_g v)_+. 
\]
Replacing $v$ by $C_g^{-1}v$ in the second equation of this pair yields 
\nonumber 
\[
Z(v_+) = (C_g^{-1}v)_- + (Z_{g^{-1}}v)_+ 
\]
whence the first equation of the pair yields 
\nonumber 
\[
Z(v_-) = (C_{g^{-1}}^{-1} v)_+ + (Z_g v)_- 
\]
on account of the identity $C_{g^{-1}}^{-1} = Z_{g^{-1}} A_g + C_g $. Conversely, reversal of this argument shows that if $Z: V_{\mathbf{C}} \longrightarrow V_{\mathbf{C}}$ is defined by these actions on $V_+$ and $V_-$ then the Gaussian $e^Z$ is the kernel of a weak metaplectic operator implementing $g$. 

\medbreak 

\begin{theorem} \label{metaplectic}
Each $g \in \mathrm{Sp}(V)$ has corresponding weak metaplectic operator $U_g: SV \longrightarrow SV'$ with kernel $u_g = e^Z$ where $Z: V_{\mathbf{C}} \longrightarrow V_{\mathbf{C}}$ is given by the rule that if $v \in V$ then 
\nonumber 
\[
Z(v_+) = (C_g^{-1}v)_- + (Z_{g^{-1}}v)_+
\]
\nonumber 
\[
Z(v_-) = (C_{g^{-1}}^{-1} v)_+ + (Z_g v)_-. 
\]
\end{theorem} 

\begin{flushright}
$\Box$
\end{flushright} 

\medbreak 

We remark that the pair of conditions on the kernel $u$ may be reformulated instructively as follows: upon addition, it follows that if $v \in V$ then 
\nonumber 
\[
\{ c(v_-) + a(v_-) \}u = c(C_g v + A_g v)_+u + a(C_g v + A_g v)_+ u
\]
whence if $\pi$ denotes the Fock representation of $V_{\mathbf{C}}$ then 
\nonumber 
\[
v \in V \Longrightarrow \pi( (g v)_+ - v_-) u = 0.
\]
Now, consider the real subspace $F_g = \{ (g v)_+ - v_- : v \in V \}$ of $V_{\mathbf{C}}$: the combined conditions on $u$ imply that if $w \in F_g$ then $a(w)u = -c(w)u$ whence with $Z$ as in Theorem 6 it follows that $a(w)\{u e^{-Z} \} = 0$ and therefore $a(i w)\{u e^{-Z} \} = 0$ by the antilinear dependence of annihilators; in view of the decomposition $V_{\mathbf{C}} = F_g \oplus i F_g$ we conclude that $ue^{-Z}$ is killed by all annihilators and is therefore a scalar. This proves that all weak metaplectic operators implementing $g$ are proportional. The specific operator $U_g$ in Theorem \ref{metaplectic} is singled out by having (generalized) vacuum expectation value $U_g \mathbf{1} (\mathbf{1}) = 1$. Indeed, let $Z: V_{\mathbf{C}} \longrightarrow V_{\mathbf{C}}$ correspond to the quadratic $\zeta \in S^2 V_{\mathbf{C}}'$ and let $Z_{g^{-1}}$ correspond to $\zeta_{g^{-1}}$. By induction along the lines of the proof for Theorem \ref{identity}, if $n \in \mathbf{N}$ and $v \in V$ then 
\nonumber
\[
2^n \:\zeta^n (v_{+}^{2n}) = (2n)! \:\Bigl(\langle v_+ \vert Z v_+ \rangle \Bigr)^n = (2n)! \:\Bigl(\langle v \vert Z_{g^{-1}} v \rangle \Bigr)^n = 2^n \:\zeta_{g^{-1}}^{n} (v^{2n})
\]
whence it follows by linearity that  
\nonumber 
\[
\psi \in S^{2n} V \Longrightarrow \zeta^n (\psi_+) = \zeta_{g^{-1}}^{n} (\psi)
\]
and therefore by summation that 
\nonumber 
\[
U_g \mathbf{1} = e^{Z_{g^{-1}}}.
\]

\medbreak 

Now, if $U_g : SV \longrightarrow SV'$ actually maps into Fock space $S[V]$ then in particular $e^{Z_{g^{-1}}} = U_g \mathbf{1} \in S[V]$ whence our knowledge of Gaussians reveals that $Z_{g^{-1}}$ is of Hilbert-Schmidt class. In the opposite direction, let us suppose that $Z_{g^{-1}}$ is Hilbert-Schmidt ; equivalently $Z_g$ is Hilbert-Schmidt, on account of the identity $Z_{g^{-1}} C_g = - C_g Z_g$. Again our knowledge of Gaussians steps in, informing us that $U_g \mathbf{1} = e^{Z_{g^{-1}}}$ lies in $S[V]$ and indeed in the domain of each creator-annihilator polynomial, whence $U_g$ maps the whole of $SV$ into $S[V]$ on account of the intertwining property. We may therefore define a Hermitian form $[ \cdot \vert \cdot]$ on $SV$ by the rule that if $\phi, \psi \in SV$ then 
\nonumber 
\[
 [\psi \vert \phi] = \langle U_g \psi \vert U_g \phi \rangle.
\]
If also $v \in V$ then the mutually adjoint nature of creators and annihilators yields 
\nonumber 
\[
\langle \{ c(C_g v) + a(A_g v) \} U_g \psi \vert U_g \phi \rangle = \langle \ U_g \psi \vert \{ a(C_g v) + c(A_g v) \} U_g \phi \rangle
\]
whence the intertwining property of $U_g$ yields 
\nonumber 
\[
[c(v) \psi \vert \phi] = [\psi \vert a(v) \phi]. 
\]
According to Lemma \ref{innerproduct} it now follows that $U_g : SV \longrightarrow S[V]$ is isometric up to the scalar 
\nonumber 
\[
[\mathbf{1} \vert \mathbf{1}] = \Vert e^{Z_{g^{-1}}} \Vert^2 = \mathrm{Det}^{\frac{1}{2}}(I - Z_{g^{-1}}^2)^{-1} = \mathrm{Det}^{\frac{1}{2}}(I - Z_{g}^2)^{-1}
\]
thus, 
\nonumber 
\[
\phi, \psi \in SV \Longrightarrow \langle U_g \psi \vert U_g \phi \rangle =  \mathrm{Det}^{\frac{1}{2}}(I - Z_{g}^2)^{-1} \:\langle \psi \vert \phi \rangle. 
\]
In similar fashion, $U_{g^{-1}}$ maps $SV$ to $S[V]$ and is isometric up to the same scalar; moreover, direct calculation reveals the relation $U_{g}^* = U_{g^{-1}}$ in the sense of Theorem \ref{adjoint}. Thus, the continuous extension $\overline{U_g} : S[V] \longrightarrow S[V]$ is (up to the very same scalar) both isometric and coisometric, hence unitary. In this way, we recover the standard Shale theorem in the following form. 

\medbreak 

\begin{theorem} \label{Shale} 
The symplectic automorphism $g \in \mathrm{Sp}(V)$ admits a unitary metaplectic operator precisely when $Z_g$ is of Hilbert-Schmidt class; a canonical choice of unitary metaplectic operator is $\mathrm{Det}^{\frac{1}{4}}(I - Z_{g}^2) \overline{U_g}$. 
\end{theorem}

\medbreak 

It is worth indicating briefly how the metaplectic kernels obtained in this paper reproduce those obtained elsewhere in the literature. As in Theorem \ref{metaplectic}, let $g \in \mathrm{Sp}(V)$ have metaplectic kernel $u_g = e^Z$ and again let $Z$ correspond to the quadratic $\zeta$. If $x, y \in V$ then inductively 
\nonumber 
\[
n \in \mathbf{N} \Longrightarrow 2^n \: \zeta^n \Bigl((x_+ + y_-)^{2n}\Bigr) = (2n)! \: \Bigl\{ \langle (x_+ + y_-) \vert Z (x_+ + y_-) \rangle \Bigr\}^n 
\]
whence by summation 
\nonumber 
\[ 
u_g \Bigl( (e^x)_+ (e^y)_- \Bigr) = u_g \Bigl( e^{x_+ + y_-} \Bigr) = \exp\Bigl\{ \frac{1}{2}\langle (x_+ + y_-) \vert Z (x_+ + y_-) \rangle \Bigr\}; 
\] 
as a direct calculation from the formula for $Z$ in Theorem \ref{metaplectic} shows that 
\nonumber 
\[ 
\langle (x_+ + y_-) \vert Z (x_+ + y_-) \rangle = 2\langle C_{g}^{-1} x \vert y \rangle + \langle x \vert Z_{g^{-1}} x \rangle + \langle Z_g y \vert y \rangle
\]
we have reproduced the standard metaplectic kernel formula. Although the expression $u_g \Bigl( e^{x_+ + y_-} \Bigr)$ appears purely formal, the simple exponentials $e^x$ and $e^y$ actually lie in the Fock space $S[V]$ and when $Z_g$ is Hilbert-Schmidt there results the authentic formula 
\nonumber 
\[ 
\langle e^x \vert \overline{U_g} e^y \rangle = \exp \frac{1}{2}\Bigl\{2\langle C_{g}^{-1} x \vert y \rangle + \langle x \vert Z_{g^{-1}} x \rangle + \langle Z_g y \vert y \rangle \Bigr\}. 
\]

\medbreak 

We close by remarking on the situation as regards implementation of antisymplectics. By definition, an antisymplecic automorphism of $V$ is a real-linear automorphism $g$ with the property that if $x,y \in V$ then $\Omega(g x,g y) = - \Omega(x,y)$; note that in this case, the antilinear part $A_g$ is invertible and $Z_g = A_{g}^{-1}C_g$ is symmetric with operator norm strictly less than unity. By a (weak) metaplectic operator implementing $g$ we shall mean a (nonzero) antilinear map $U: SV \longrightarrow SV'$ that intertwines in the sense that if $v \in V$ then 
\nonumber 
\[ 
U c(v) = \{ a(C_g v) + c(A_g v) \} U 
\] 
\nonumber 
\[ 
U a(v) = \{ a(A_g v) + c(C_g v) \} U. 
\] 
Making use of the kernel calculus in Theorem \ref{antikernel} it follows that the antisymplectic $g$ always admits the weak metaplectic operator $U_g$ with kernel $u_g = e^Z \in S(V \oplus V)'$ where $Z: V \oplus V \longrightarrow V \oplus V$ is given by the rule that if $v \in V$ then 
\nonumber 
\[ 
Z(v_1) = (Z_g v)_1 + (A_{g^{-1}}^{-1} v)_2 
\] 
\nonumber 
\[ 
Z(v_2) = (A_{g}^{-1} v)_1 + (Z_{g^{-1}} v)_2 .
\] 
This weak metaplectic operator $U_g$ is unique up to scale; it may be rescaled so as to be antiunitary precisely when $C_g$ is of Hilbert-Schmidt class. 

\medbreak 

In short: just as (Shale) a symplectic is implemented by a unitary operator on Fock space precisely when its commutator with $i$ is Hilbert-Schmidt, so (anti-Shale) an antisymplectic is implemented by an antiunitary on Fock space precisely when its anticommutator with $i$ is Hilbert-Schmidt.

\end{document}